\documentclass[12pt]{amsart}

\newcommand{\R}{{\mathbb R}}

\newcommand{\Z}{{\mathbb Z}}

\renewcommand{\phi}{\varphi}

\newtheorem{theo}{{\sc Theorem}}

\newtheorem{lem}[theo]{{\sc Lemma}}
\newtheorem{prop}[theo]{{\sc Proposition}}

\newenvironment{defn}{\medskip\noindent{\it Definition:\/} }{\medskip}

\title[$L^p$ norms of eigenfunctions in the completely integrable case]{$L^p$  norms of eigenfunctions in
the completely integrable case}

\author{John A. Toth and Steve Zelditch}
\address{Department of Mathematics and Statistics, McGill University,
Montreal, CANADA, H3A-2K6}
\address{Department of Mathematics, Johns Hopkins University, Baltimore, MD
21218, USA}

\thanks{\\
 Research partially supported by an Alfred P. Sloan Research Fellowship and  NSERC grant \#OGP0170280\\ Research
partially supported by  NSF grant \#DMS-0071358}

\date{\today}

\begin{document}

\maketitle

\addtolength{\baselineskip}{1pt}

\begin{abstract} The eigenfunctions $e^{i \langle \lambda, x
\rangle}$ of the Laplacian on a flat torus have uniformly bounded
$L^p$ norms.   In this article, we prove that for every other
quantum integrable Laplacian, the $L^p$ norms of the joint
eigenfuntions blow up at least at the rate  $\| \phi_{k}
\|_{L^{p}} \geq C(\epsilon) \lambda_{k}^{ \frac{ p - 2}{ 4p } -
\epsilon }$ when $p >2$. This gives a quantitative refinement of our recent
result \cite{TZ1} that some sequence of eigenfunctions must blow
up in $L^p$ unless $(M, g)$ is flat. The better result in this
paper is based on mass estimates of eigenfunctions near singular
leaves of the Liouville foliation.

\end{abstract}

\setcounter{page}{1} \setcounter{section}{-1}
\section{Introduction}

This paper, a companion to \cite{TZ1},  is concerned with the
growth rate of the $L^{p}$-norms of  $L^2$-normalized
$\Delta$-eigenfunctions
$$\Delta \phi_j = \lambda_j^2 \phi_j, \;\;\;\; \langle \phi_j,
\phi_k \rangle = \delta_{jk}$$
 on compact Riemannian manifolds $(M, g)$ with completely integrable geodesic
 flow $G^t$ on $S^*M$.  The motivating problem is to relate  sizes of
eigenfunctions  to dynamical properties of its geodesic flow $G^t$
on $S^*M$. In general this is an intractable problem, but much can
be understood by studying it in the framework of integrable
systems.  To be  precise, we assume that   $\Delta$ is
  {\it quantum completely integrable} or QCI in the sense
that there exist $P_1, \dots, P_n \in \Psi^1(M)$ ($n$ = dim $M$)
 satisfying
$[P_i, P_j]=0$ and such that their symbols $(p_1, \dots, p_n)$
satisfy the standard assumption that  $dp_1\wedge dp_2 \wedge
\dots \wedge dp_n \not= 0$ on a dense open set $\Omega \subset
T^*M - 0$ whose complement is contained in a hypersurface. Since
$\{p_i, p_j\} = 0$, the $p_1, \dots, p_n$ generate a Hamiltonian
$\R^n$-action given by the joint Hamiltonian flow
$\Phi_{t}(x,\xi):= \exp(t_{1} X_{p_{1}}) \circ ... \circ
\exp(t_{1} X_{p_{n}}) (x,\xi )$, where for $j=1,...,n$,
$X_{p_{j}}$ denotes the Hamiltonian vector field of $p_{j}$. The
associated moment map ${\mathcal P}: T^{*}M \rightarrow {\Bbb
R}^{n}$ is  ${\mathcal P}= (p_{1},...,p_{n})$ and the orbits of
this action foliate $T^*M - 0$ and since $S^*M$ is preserved they
also foliate $S^*M$.  See \cite{TZ1, TZ2} for a list of well-known
examples. We will refer to this foliation as the {\it Liouville
foliation}. Throughout, we will
 assume that the leaves of this foliation satisfy  Eliasson's
non-degeneracy hypotheses (see \cite{El} and Definition \ref{END}).

The main result of \cite{TZ1} was that the $L^{\infty}$-norms of
the  $L^2$-normalized joint eigenfunctions $\{\phi_{\lambda}\}$ of
$(P_1, \dots, P_n)$ are unbounded unless  $(M, g)$ is flat.  In
this paper, we investigate the blow-up rates of $L^{p}$-norms of
sequences $||\phi_{\lambda}||_{L^{p}} $, under the Eliasson
non-degeneracy assumption.

\begin{theo} \label{LP} Suppose that $(M, g)$ is a compact  Riemannian manifold
with completely integrable geodesic flow  satisfying Eliasson's
non-degeneracy condition. Then, unless  $(M, g)$ is a flat torus,
  there exists  for every $\epsilon >0$, a sequence of
eigenfunctions satisfying:
$$\left\{ \begin{array}{l}   \| \phi_{k} \|_{L^{\infty}} \geq C(\epsilon) \lambda_{k}^{\frac{1}{4} -
\epsilon}. \\ \\
 \| \phi_{k} \|_{L^{p}} \geq C(\epsilon) \lambda_{k}^{ \frac{ p - 2}{ 4p } - \epsilon }, \,\,\, 2 < p < \infty. \end{array} \right.$$
   \end{theo}

This result is sharp in the setting of all completely integrable
systems. It is based on the existence of a codimension one
singular leaf of the Liouville foliation. Here, a singular leaf is
an orbit of the $\R^n$ action with dimension $< n$, see Definition
(\ref{SING}).  When codimension $\ell$ leaves occur, the estimate
becomes $\| \phi_{k} \|_{L^{\infty}} \geq C(\epsilon)
\lambda_{k}^{\frac{\ell }{4} - \epsilon}$

 Existence of a singular leaf is given in the following:

\begin{lem} \label{FLAT} Suppose that $(M, g)$ is a compact Riemannian manifold with completely integrable geodesic flow.
Then the  Liouville foliation of $(M, g)$ contains a  singular
leaf unless $(M, g)$ is a flat torus.
\end{lem}

The proof only shows that some singular leaf occurs, but does not
determine its codimension. As will be discussed in \S
\ref{GEOMETRY}, it is quite plausible that codimension $n - 1$
singular leaves  often  occur for Hamiltonian $\R^n$ actions on
cotangent bundles. Such leaves correspond to closed geodesics
which are invariant under the $\R^n$ action. In such cases, the
estimate improves to $ \| \phi_{k} \|_{L^{\infty}} \geq C(\epsilon)
\lambda_{k}^{\frac{n - 1}{4} - \epsilon}.$ This estimate agrees
with a recent result of Donnelly \cite{D} on blow-up rates of
$L^{\infty}$ norms of eigenfunctions on compact $(M, g)$ with
isometric $S^1$ actions. He finds sequences of joint
eigenfunctions which blow up at the rate $\lambda_{k}^{\frac{n -
1}{4} }$; they are precisely the same kind of eigenfunctions which
are associated to singular leaves, as we now describe.

 The analytic ingredient of our proof  gives an estimate on the growth rate of
modes and quasimodes associated to a singular leaf of the
Liouville foliation. It is here we assume that the singular leaf
is Eliasson non-degenerate.

\begin{lem} \label{BLOWUP}   Let ${\mathcal P}^{-1}(c)$ with $c= (c_{1},c_{2},...,c_{n}) \in \R^{n}$ be a singular level
of the Lagrangian fibration associated with the quantum completely
integrable system $P_{1}= \sqrt{\Delta},...,P_{n}$. Suppose
${\mathcal P}^{-1}(c)$ contains a single compact, Eliasson non-degenerate
orbit $\Lambda:= {\Bbb R}^{n} \cdot (v_{0})$ of dimension $\ell
<n$ and that each connected component of $\pi( \Lambda ) \subset
M$ is an embedded $\ell$-dimensional submanifold. Then, there
exists a sequence  $\phi_{k} $ of   $L^{2}$-normalized joint
eigenfunctions of $P_{1},...,P_{n}$ with $P_{j} \phi_{k} = E_{jk}
\phi_{k} $, \,  $E_{jk} = \lambda_{k} c_{j} + {\mathcal O}(1); \,
j=1,...,n$ \, such that for any $\epsilon >0$,

$$  (i) \,\, \| \phi_{k} \|_{L^{\infty}} \geq C(\epsilon) \lambda_{k}^{\frac{n-\ell}{4} -
\epsilon}.$$
We also have that
$$ (ii) \,\, \| \phi_{k} \|_{L^{p}} \geq C(\epsilon) \lambda_{k}^{ \frac{ (n- \ell ) (p-2)}{ 4 p }  - \epsilon }, \,\,\, 2 < p < \infty.$$
\end{lem}

For simplicity we have stated Lemma (\ref{BLOWUP}) in the case where the
level ${\mathcal P}^{-1}(c)$ contains a single Eliasson
non-degenerate orbit, but the result extends in a straightforward
fashion to the case of multiple singular orbits.

Some final words to close the introduction. First, in \cite{TZ2}
we give a different approach to eigenfunction blow-up in the
completely integrable case which in some ways is closer in spirit
to the approach of this paper.  Namely, we relate norms of  modes
to norms of ``quasi-modes,'' i.e.  approximate eigenfunctions
associated to Bohr-Sommerfeld leaves of the Liouville foliation.
In that paper, we also give a number of detailed examples of
quantum completely integrable systems such as Liouville tori,
surfaces of revolution, ellipsoids, tops and so on.

Second, we would like to emphasize that results on completely
integrable systems are relevant to quite  general classes of
Riemannian manifolds. This is because one can approximate any
geodesic flow by a completely integrable flow (its Birkhoff normal
form) near a closed orbit or invariant torus. On the quantum
level, one may approximate the Laplacian by its quantum Birkhoff
normal form.  The quasi-mode analysis of this paper is therefore
of some relevance to general Riemannian manifolds. What is
 special about quantum integrable systems is the particularly
 close relations between modes and quasi-modes, which allows one
 to draw strong conclusions about the blow-up of sequences of
 modes from that of the corresponding quasi-modes.  It is quite
 reasonable to suppose that the same kind of phenomenon occurs
 much more generally. In particular, we would conjecture that any
 compact $(M, g)$ with a stable elliptic orbit has a sequence of
 eigenfunctions whose $L^{\infty}$ norms blow up at the rate
 $\lambda_{k}^{\frac{n - 1}{4} }.$  There surely exist quasi-modes
 with this property, and the difficulty is to show that this
 implies the existence of modes with this blow-up rate. We hope to
 pursue this direction in the future.

\section{\label{GEOMETRY} Geometry of completely integrable systems}

This section is devoted to the geometric aspects of our problem.
 We first  prove Lemma (\ref{FLAT}) and then give
background on Eliasson non-degeneracy of singular orbits. We begin
with  some preliminary background on completely integrable
systems, in part repeated from  \cite{TZ1}.

  As mentioned in the introduction, the
moment map of a completely integrable system is defind by
\begin{equation} \label{MM} {\mathcal P} = (p_1, \dots, p_n): T^*M \rightarrow
B \subset\R^n. \end{equation}   The Hamiltonians $p_j$ generate
the
 $\R^n$-actin
$$ \Phi_t = \exp t_1 \Xi_{p_1} \circ \exp t_2 \Xi_{p_2} \dots \circ \exp t_n
\Xi_{p_n}.$$ We  denote  $\Phi_t$-orbits by   $\R^n \cdot (x, \xi)
$.

By the Liouville-Arnold theorem, the orbits of the joint flow
$\Phi_{t}$ are diffeomorphic to $\R^k \times T^m$ for some $(k,m),
k + m \leq n.$ By the properness assumption on  ${\mathcal P}$,
each connected component of a regular level is a Lagrangean torus,
i.e.
\begin{equation} \label{CI1}
{\mathcal P}^{-1}(b) = \Lambda^{(1)}(b) \cup \cdot \cdot \cdot
\cup
 \Lambda^{(m_{cl})}(b) , \;\;\;(b \in B_{reg})
\end{equation}
\noindent where each $\Lambda^{(l)}(b) \simeq T^n$ is a Lagrangian
torus.

 We are particularly interested in singular orbits. To clarify
the notation and terminology, we repeat the definition from
\cite{TZ1}:

\begin{defn}\label{SING}
We say that:

\begin{itemize}

\item  $b \in B_{sing}$ if ${\mathcal P}^{-1}(b)$ is a singular level of the moment map, i.e. if
there exists a point $(x, \xi) \in {\mathcal P}^{-1}(b)$  with
$dp_{1} \wedge \cdot \cdot \cdot \wedge dp_{n}(x,\xi) = 0$. Such a
point $(x, \xi)$ is called a singular point of ${\mathcal P}$.

\item a connected component of ${\mathcal P}^{-1}(b)$ ($b \in B_{sing}$)is a singular component if it  contains a singular point ;

\item an orbit $\R^n \cdot (x, \xi) $ of $\Phi_{t}$ is  singular if it is non-Lagrangean, i.e.   has dimension $<n$;

\item $b \in B_{reg}$ and that ${\mathcal P}^{-1}(b)$ is a regular level if all points $(x, \xi) \in {\mathcal P}^{-1}(b)$ are regular, i.e. if $dp_{1} \wedge \cdot \cdot \cdot \wedge dp_{n}(x, \xi) \not= 0$.

\item   a component of  ${\mathcal P}^{-1}(b)$ ( $b \in B_{sing} \cup B_{reg}$) is regular if it contains no singular points.

\end{itemize}

\end{defn}

 When  $b \in B_{sing}$ we first
decompose
\begin{equation} \label{SL} {\mathcal P}^{-1}(b) = \cup_{j = 1}^{r} \Gamma_{sing}^{(j)}(b) \end{equation}
the singular level into connected components
$\Gamma_{sing}^{(j)}(b)$ and then decompose
\begin{equation} \label{ORB} \Gamma_{sing}^{(j)}(b) = \cup_{k = 1}^{p} \R^n \cdot (x_k, \xi_k) \end{equation}
each component  into orbits. Both decompositions can take a
variety of forms.  The regular components $\Gamma_{sing}^{(j)}(b)$
must be Lagrangean tori by the properness assumption. A singular
component consists of finitely many orbits by the Eliasson non-degeneracy
assumption (see below). The orbit $\R^n \cdot (x, \xi)$ of a
singular point is necessarily singular, hence has the form $\R^k
\times T^m$ for some $(k,m)$ with  $k + m < n.$ Regular points may
also occur on a singular component, whose orbits are Lagrangean
and can take any one of the forms $\R^k \times T^m$ for some
$(k,m)$ with  $k + m = n.$

\subsection{Proof of Lemma (\ref{FLAT})}

In this section, we prove Lemma (\ref{FLAT}):  the Liouville
foliation of  metrics with completely integrable geodesic flows
always contain singular leaves unless $(M, g)$ is a flat torus.

\begin{proof} The hypothesis that the Liouville foliation is
non-singular has two immediate  geometric consequences:
\begin{itemize}

\item (i) By Mane's theorem \cite{M} , $(M, g)$ is a manifold without
conjugate points;

\item (ii)  The (homogeneous) moment map ${\mathcal P}: T^*M-0 \to \R^n-0$ is
a torus fibration by $T^n$.

\end{itemize}

Statement (ii) follows from the Liouville-Arnold theorem. On a
non-singular leaf, we must have
 $dp_1 \wedge \cdots \wedge dp_n \not= 0$. Since this holds everywhere,  ${\mathcal P}$ is a
 submersion; and since it is proper, it is a fibration. The fiber
 must be  $T^n$, again  by the Liouville-Arnold theorem. Since ${\mathcal
 P}$ is homogeneous, the image ${\mathcal P}(T^*M - 0) = \R_+
 \cdot {\mathcal P}(S^*M)$. Since ${\mathcal P}$ is a submersion,
 the image is a smooth submanifold of $S^{n-1}$ hence must be all of
 $S^{n-1}$.

By (i) it follows that $M = \tilde{M}/\Gamma$ where $(\tilde{M},
\tilde{g})$ is the universal Riemannian cover of $(M, g)$
(diffeomorphic to $\R^n$), and where $\Gamma \cong \pi_1(M)$ is
the group of covering transformations.

We claim that (ii) implies $\pi_1(M) = \Z^n$. Indeed, $T^*M-0$ is
a double fibration (with indicated fibers):
$$\begin{array}{lllll} & & T^*M-0 & & \\& &  & & \\
& \pi \; \swarrow \; (\R^n-0)   & & (T^n) \; \searrow \; {\mathcal P} & \\ & & &  & \\
M & & & & \R^n-0 \end{array}. $$

By the homotopy sequence of a  fibration $\pi: E \to B$,
$$\cdots \pi_q(F) \rightarrow \pi_q(E) \to
\pi_q(B) \to \pi_{q-1}(F) \cdots \to \pi_0(E) \to \pi_0(B) \to 0
$$
and using that $\pi_2(T^n) = {\bf 1}$ and that $\pi_2(M) = {\bf
1}$ by (i), we obtain
\begin{equation} \begin{array}{l} {\bf 1}  \to  \pi_1(S^{n-1}) \to \pi_1(S^*M) \to
\pi_1(M) \to  \pi_{0}(S^{n-1}) \to \pi_0(S^*M) \to \pi_0(M) \to {\bf 1} \\ \\
{\bf 1}  \to \pi_2(S^*M) \to \pi_2(S^{n-1}) \to \pi_1(T^n)
\rightarrow \pi_1(S^*M) \to \pi_1(S^{n-1}) \to \pi_{0}(T^n)  \to
\pi_0(S^*M) \to \pi_0(S^{n-1} \to  {\bf 1}
\end{array} \end{equation}

Since $\pi_2(N) = \pi_2(\tilde{N})$ (where $\tilde{N}$) is the
universal cover, and since $\tilde{S^*M} = \tilde{S^{n-1}} \times
\R^n$,  we have $\pi_2(S^*M) = \pi_2(S^{n-1})$. From its
 definition,  we see that the homomorphism $ \pi_2(S^*M) \to \pi_2(S^{n-1})$
 is an isomorphism, hence
 the second sequence simplifies to
\begin{equation} {\bf 1} \to \pi_1(T^n)
\rightarrow \pi_1(S^*M) \to \pi_1(S^{n-1}) \to \pi_{0}(T^n)  \to
\pi_0(S^*M) \to \pi_0(S^{n-1} \to  {\bf 1}
\end{equation}

Let us first assume that $n \geq 3$. Then $\pi_1(S^*M) =
 \pi_1(M)$ and $\pi_1(S^{n-1}) = {\bf 1}$, so  we get
 $${\bf 1} \to   \pi_1(T^n) \rightarrow \pi_1(M) \to
 {\bf 1}, $$
i.e. $\pi_1(T^n) \cong
 \pi_1(M).$

We now consider dimension $2$. Actually, it follows by a classic
result of Kozlov \cite{K} that the only surfaces that can possibly
have a completely integrable geodesic flow (even with
singularities) are $M = S^2, T^2$. We can easily disqualify $M =
S^2$ since $S^*S^2 \equiv \R P^3 $ does not fiber over $S^1$.

It follows then that $M$ is diffeomorphic to $\R^n/\Z^n$, i.e. $M$
is a torus. Since $g$ has no conjugate points, the proof is
concluded by Burago-Ivanov's theorem \cite{BI}.

\end{proof}

\subsubsection{Remarks on the geometry} We would like to add two
remarks on the geometry. We will not use them in the proof of our
results, and therefore only sketch the issues. Our second remark
could lead to a substantial improvement of our main result.

\noindent{\bf Remark (i): On flat manifolds.}

In our previous article \cite{TZ1}, we showed that eigenfunction
blow-up occurs in the integrable setting unless $(M, g)$ is flat.
Here, we are showing that it occurs unless $(M, g)$ is a flat
torus.   We briefly explain why Liouville foliations of other flat
manifolds must have singular leaves.

The Riemannian connection of $(M, g)$ gives a distribution of
horizontal planes in the unit sphere bundle, i.e. $T(S^*M) = H
\oplus V$ where $V$ is the tangent bundle to the fibers of $\pi:
S^*M \to M$. If $g$ is a flat metric, then $H$ is involutive and
hence $S^*M$ is foliated by smooth $n$-manifolds $L$ (which we
will call horizontal leaves). $(M, g)$ is also the Riemannian
quotient of a torus $T^n$ by a finite group $G$ of isometries, and
the derivative $p_*$ of the  projection $p : T^n \to M$ takes the
horizontal distribution of $T^n$ to that of $M$. Since the
horizontal leaves of $T^n$ are compact Lagrangean tori, it follows
that the horizontal leaves of $(M, g)$ are compact Lagrangean
submanifolds of $T^*M - 0.$  Thus, the geodesic flow of $(M, g)$
leaves invariant a non-singular Lagrangean foliation.

However, the horizontal leaves are not orbits of a Hamiltonian
$\R^n$ action in general. (We thank W. Goldman for several helpful
discussions on this point).  In fact, the leaves have the form
$L(G)\backslash \R^n/\Z^n$ where $L(G)$ is the linear holonomy
group (a finite group). So the leaves are quotients of tori, but
$L(G)$ is generally not a subgroup of $T^n$ so the leaves are
finite quotients of tori. If $L(G)$ acts freely on $T^*M - 0$,
then the leaves are tori, but there is no $\R^n$ action generating
them.

The natural `Hamiltonians' on $T^*M - 0$ belong to the algebra
${\mathcal I}^G$ where ${\mathcal I} = <I_1, \dots, I_n>$ is the
algebra generated by the action variables $I_j (x, \xi) = \xi_j$
on the covering torus, and where ${\mathcal I}^G$ denotes the
$G$-invariant elements. Since we require that the $\R^n$ action be
homogeneous, we only consider functions $q_j(I_1, \dots, I_n)$
which are homogeneous of degree one relative to the natural $\R_+$
action on $T^*M-0.$ A completely integrable system on $T^*M - 0$
is given by a choice of $n$ generators of this algebra. The $\R^n$
action they generate will leave invariant the horizontal
foliation, but the leaves of that foliation are not in general the
orbits of the action.  What the proof above shows is that any such
action must have some singular orbits, of lower dimension than the
leaf it is contained in.
\medskip

\noindent{\bf Remark (ii): On existence of high codimension
singular leaves. }

The simple proof above on existence of singular leaves does not
give any information on the codimension of such leaves. As
mentioned in the introduction, our lower bounds
  $||\phi_k||_{L^{\infty}} \geq C(\epsilon)
\lambda_k^{\frac{q}{4} - \epsilon}$ improve as this codimension
$q$ increases. Can a better argument prove the existence of
codimension $n - 1$ leaves in generic cases?

Here is heuristic existence argument: the image ${\mathcal C}$ of
$T^*M - 0$ under the moment map ${\mathcal P}$ of a homogeneous
$\R^n$-action is a convex polyhedral cone in $\R^n$. In Lemma
(\ref{FLAT}), we observed  that ${\mathcal P}$ cannot be a torus
bundle over $\R^n- 0$ unless $(M, g)$ is a flat torus. In all
other cases, ${\mathcal C}$ will have a boundary over which the
singular orbits `fiber'.

The geometric problem is to determine the geometry of ${\mathcal
C}$. In the case of Hamiltonian torus actions, the moment
polyhedra have been studied by Lerman et al. \cite{L, L2}. They do
not seem to have been studied for general Hamiltonian $\R^n$
actions. However, it is natural to conjecture that as long as
${\mathcal C}$ has a non-empty boundary, there will  boundary
faces of each codimension and in particular there will be
one-dimensional edges. If this is true, then our lower bound
improves to $ >> \lambda_k^{\frac{n-1}{4} - \epsilon}$.

\subsection{\label{END} Eliasson non-degeneracy}

Let $c:= (c_{1},...,c_{n}) \in B_{sing}$ be a singular  value of
the moment map, ${\mathcal P}$. Suppose that:
\begin{equation} \label{RANK}
{\mbox rank}\, (dp_{1},...,dp_{n}) = {\mbox rank}\,
(dp_{1},...,dp_{k}) = k < n
\end{equation}
\noindent at some point $ v_{0} \in {\mathcal P}^{-1}(c)$. Denote
the orbit through $v_{0}$ by ${\Bbb R}^{n} \cdot (v_{0}) := \{
\exp t_{1} \Xi_{p_{1}} \circ \cdot \cdot \cdot \circ \exp t_{k}
\Xi_{p_{k}} (v_{0}); \, \, t = (t_{1},...,t_{k}) \in {\Bbb R}^{n}
\}$, which we henceforth assume is compact. Our aim is to define
the notion of Eliasson non-degeneracy of such a singular orbit;
this is a straightforward extension of the notion defined in
\cite{El} for fixed points of symplectic maps.

By the Liouville-Arnold theorem, the orbits of the joint flow
$\Phi_{t}$ are diffeomorphic to $\R^k \times T^m$ for some $(k,m),
k + m \leq n.$ Let us first consider regular levels where $m = n$.
By the properness assumption on ${\mathcal P}$, a regular level
has the form
\begin{equation} \label{CI1}
{\mathcal P}^{-1}(b) = \Lambda^{(1)}(b) \cup \cdot \cdot \cdot
\cup
 \Lambda^{(m_{cl})}(b) , \;\;\;(b \in B_{reg})
\end{equation}
\noindent where each $\Lambda^{(l)}(b) \simeq T^n$ is an
$n$-dimensional Lagrangian torus.  The  classical (or geometric)
multiplicity function   $m_{cl} (b) = \# {\mathcal P}^{-1} (b)$,
i.e. the number of  orbits on the level set ${\mathcal
P}^{-1}(b)$, is constant on    connected components of $B_{reg}$
and the moment map ${\mathcal P}$ is a fibration over each component
with fiber (\ref{CI1}).
 In  sufficiently small neighbourhoods $\Omega^{(l)}(b)$ of each component torus,
$\Lambda^{(l)}(b)$, the Liouville-Arnold theorem also gives the
existence of local action-angle variables
$(I^{(l)}_{1},...,I^{(l)}_{n},
\theta^{(l)}_{1},...,\theta^{(l)}_{n})$ in terms of which the
joint flow of $\Xi_{p_{1}},...,\Xi_{p_{n}}$ is linearized [AM].
For convenience, we henceforth normalize the action variables
$I^{(l)}_{1},...,I^{(l)}_{n}$ so that $I^{(l)}_{j} = 0; \,
j=1,...,n$ on the torus $\Lambda^{(l)}(b)$.

Now let us consider singular levels of rank $k$.  We first observe
that $dp_{1},...,dp_{k}$ (in the notation of (\ref{RANK})) are
linearly independent everywhere on ${\Bbb R}^{n} \cdot (v_{0})$.
Indeed,  the Liouville-Arnold theorem in the singular case [AM]
states that  there exists  local canonical transformation
$$ \psi = \psi(I,\theta,x,y): {\Bbb R}^{2n} \rightarrow
T^{*}M-0,$$ where $$ I=(I_{1},...,I_{k}),
\theta=(\theta_{1},...,\theta_{k}) \in {\Bbb R}^{k},\;
x=(x_{1},...,x_{n-k}), y=(y_{1},...,y_{n-k}) \in {\Bbb R}^{n-k}m
$$ \noindent defined in an invariant neighbourhood of ${\Bbb R}^{n} \cdot (v_{0})$ such that
\begin{equation}
p_{i} \circ \psi = I_{i} \,\,\,(i=1,...,k),
\end{equation}
and such that the symplectic form $\omega$ on $T^*M$ takes the
form
\begin{equation} \label{E1}
\psi^{*}\omega = \sum_{j=1}^{k} dI_{j} \wedge d\theta_{j} +
\sum_{j=1}^{n-k} dx_{j} \wedge dy_{j}.
\end{equation}
\noindent As a consequence of the normal form (\ref{E1}), it
follows that there exist constants $c_{ij}$ with  $i = k+1,...,n$
and $j=1,...,k$,  such that at each point of the orbit, ${\Bbb
R}^{n} \cdot (v_{0})$,
\begin{equation} \label{span}
dp_{i} = \sum_{j=1}^{k} c_{ij} dp_{j}.
\end{equation}

\noindent Since  $dp_{1},...,dp_{k}$ are linearly independent in a
sufficiently neighbourhood $U$ of $v_{0}\in {\mathcal P}^{-1}(c)$,
the action of the flows corresponding to the Hamilton vector
fields, $\Xi_{p_{1}},...,\Xi_{p_{k}}$ generates a symplectic
${\Bbb R}^{k}$ action on ${\mathcal P}^{-1}(c_{0}) \cap U$. An
application of Marsden-Weinstein reduction yields an open
$2(n-k)$-dimensional symplectic manifold:
\begin{equation} \label{E2}
\Sigma_{k} := {\mathcal P}^{-1}(c_{0}) \cap U / {\Bbb R}^{k},
\end{equation}

\noindent with the induced symplectic form, $\sigma$. We will
denote the canonical projection map by:
$$ \pi_{k} : {\mathcal P}^{-1}(c_{0}) \cap U \longrightarrow \Sigma_{k}.$$
\noindent Since $\{p_{i},p_{j} \} =0$ for all $i,j =1,...,n$, it
follows that $p_{k+1},...,p_{n}$ induce $C^{\infty}$ functions on
$\Sigma_{k}$, which we will, with some abuse of notation, continue
to write as $p_{k+1},...,p_{n}$. From (\ref{span}), it follows
that
$$dp_{i}(\pi_{k}(v_{0})) = 0;\,\,\,i=k+1,...,n.$$
\noindent Here, we denote the single point $\pi_{k}({\Bbb R}^{n}
\cdot(v_{0}))$ by $\pi_{k}(v_{0})$.

 To describe the Eliasson construction, let $(M^{2n},\sigma)$ be any symplectic manifold and ${\mathcal C}$,  the Lie
algebra of Poisson commuting functions, $f_{i} \in C^{\infty}(M)$,
with the property that:
$$ f_{i}(m) = df_{i}(m) = 0; \,\,\,i=1,...n,$$
\noindent for a fixed $m \in M$. Let ${\mathcal Q}(2n)$ denote the
Lie algebra of quadratic forms and consider the Lie algebra
homomorphism, ${\mathcal H}: {\mathcal C} \rightarrow {\mathcal
Q}(2n)$ given by
$$ {\mathcal H}(f) := d^{2} f(m).$$

\noindent The result of Eliasson \cite{El} says that, if
${\mathcal H}(f_{1},...,f_{n})$ is a Cartan subalgebra of
${\mathcal Q}(2n)$, there exists a locally-defined canonical
mapping, $\kappa : U \rightarrow U_{0}$, from a neighbourhood,
$U$,  of $m \in M$ to a neighbourhood, $U_{0}$, of $0\in {\Bbb
R}^{n}$, with the property that:
\begin{equation} \label{E3}
\forall i,j \,\,\{ f_{i} \circ \kappa^{-1}, I^{e}_{j} \}= \{ f_{i}
\circ \kappa^{-1}, I^{h}_{j} \} = \{ f_{i} \circ \kappa^{-1}, I^{ch}_{j} \}
= 0.
\end{equation}
\noindent Here $I^{h}_{1},...,I^{h}_{H}, I^{ch}_{H+1},...,
I^{ch}_{H+L+1},I^{e}_{H+L+2},...,I^{e}_{2n}$ denote the standard
basis of the Cartan subalgebra \linebreak ${\mathcal
H}(f_{1},...,f_{n})$, where $I^{h}_{i}= x_{i} \xi_{i}$ in the case
of a hyperbolic summand, $I^{e}_{i} = x_{i}^{2} + \xi_{i}^{2}$ for
an elliptic summand, and $I^{ch}_{i} = x_{i} \xi_{i+1}-
x_{i+1}\xi_{i} + \sqrt{-1}( x_{i}\xi_{i} + x_{i+1} \xi_{i+1})$ in
the complex-hyperbolic case. By making a second-order Taylor
expansion about $I^{e}=I^{ch}=I^{h}=0$, it follows from (\ref{E3})
that $\forall i=1,...,n$, there locally exist $F_{ij} \in
C^{\infty}(U_{0})$ with $\{ I^{e}_{i}, F_{ij} \}=\{ I^{h}_{i},
F_{ij} \}=\{ I^{ch}_{i}, F_{ij} \}=0$ such that:

\begin{equation} \label{E4}
f_{i} \circ \kappa^{-1} = \sum_{j=1}^{H} F_{ij} \cdot I^{h}_{j} +
\sum_{j=H+1}^{H+L+1} F_{ij} \cdot I^{ch}_{j} + \sum_{j=H+L+1}^{n}
F_{ij} \cdot I^{e}_{j}.
\end{equation}
\noindent The non-degeneracy of ${\mathcal H}(f_{1},...,f_{n})$
implies that
$$ (F_{ij})(0) \in Gl(n;\R).$$

\noindent There is also a parameter-dependent version of this
result [El] that is valid in a neighbourhood of the orbit, ${\Bbb
R}^{n} \cdot(v_{0})$, and combines the normal form in (\ref{E4})
with that in (\ref{E1}). More precisely, if ${\mathcal
H}(dp_{k+1},...,dp_{n})$ is a Cartan subalgebra, then there exists
a neighbourhood, $\Omega$, of the orbit, ${\Bbb R}^{n}
\cdot(v_{0})$, and a canonical map $\kappa: \Omega \rightarrow
{\Bbb T}^{k} \times D$ with the property that, for all
$i=1,...,n$,
\begin{equation} \label{E5}
p_{i} \circ \kappa^{-1} =  \sum_{j=1}^{H} F_{ij} \cdot I^{h}_{j} +
\sum_{j=H+1}^{H+L+1} F_{ij} \cdot I^{ch}_{j} +
\sum_{j=H+L+1}^{n-k} F_{ij} \cdot I^{e}_{j} + \sum_{j=n-k+1}^{n}
F_{ij} \cdot I_{n+1-j}.
\end{equation}
\noindent Here  $I^{h}:=(I^{h}_{1},...,I^{h}_{H}), I^{ch}:=
(I^{ch}_{H+1},..., I^{ch}_{H+L+1})$ and $I^{e}:=(
I^{e}_{H+L+2},...,I^{e}_{n-k})$ denote the standard generators of
the Cartan algebra on $\Sigma_{k}$, $I:=(I_{1},...,I_{k})$ the
regular action coordinates coming from the normal form in
(\ref{E1}) and $D \subset {\Bbb R}^{n-k}$ a small ball containing
$0 \in {\Bbb R}^{n-k}$.  The $F_{ij}$ Poisson-commute with all the
action functions.

\begin{defn} \label{END}
(i) We say that the rank $
 \ell < n$ orbit, ${\Bbb R}^{n} \cdot(v_{0})$,
is {\em Eliasson non-degenerate} provided it is compact and the
algebra ${\mathcal H}(p_{\ell +1},...,p_{n})$ is Cartan.

(ii) The rank $\ell <n$ orbit ${\Bbb R}^{n} \cdot(v_{0})$ is said to
be  {\em strongly Eliasson non-degenerate} if it is Eliasson
non-degenerate and in addition, $\pi({\Bbb R}^{n} \cdot(v_{0}) )
\subset M$ is an embedded $\ell $-dimensional submanifold of $M$.
\end{defn}

\section{Quantum integrable systems and Birkhoff normal forms}

Our purpose in this section is to construct a microlocal Birkhoff
normal form for a QCI (quantum completely integrable) system near
a singular orbit. We first recall the definition of a QCI system.

\subsection{Quantum integrable systems}

Quantum completely integrable Hamiltonians (e.g. Laplacians) or
those which commute with a maximal family of observables. We now
make precise the kind of observables we will need to use.

  Given an open $U \subset {\Bbb R}^{n}$, we say that
$a(x,\xi;\hbar) \in C^{\infty}(U \times {\Bbb R}^{n})$ is  in the
symbol class $S^{m,k}(U \times {\Bbb R}^{n})$, provided
$$ |\partial_{x}^{\alpha} \partial_{\xi}^{\beta} a(x,\xi;\hbar)| \leq
C_{\alpha \beta} \hbar^{-m} (1+|\xi|)^{k-|\beta|}.$$ \noindent We
say that $ a \in S^{m,k}_{cl}(U \times {\Bbb R}^{n})$ provided
there exists an asymptotic expansion:
$$ a(x,\xi;\hbar) \sim \hbar^{-m} \sum_{j=0}^{\infty} a_{j}(x,\xi)
\hbar^{j},$$ \noindent valid for $|\xi| \geq \frac{1}{C} >0$  with
$a_{j}(x,\xi) \in S^{0,k-j}(U \times {\Bbb R}^{n})$ on this set.
We denote the associated $\hbar$ Kohn-Nirenberg quantization by
$Op_{\hbar}(a)$, where this operator has Schwartz kernel given
locally by the formula:
$$Op_{\hbar}(a)(x,y) = (2\pi \hbar)^{-n} \int_{{\Bbb R}^{n}}
e^{i(x-y)\xi/\hbar} \,a(x,\xi;\hbar) \,d\xi.$$ By using a
partition of unity, one constructs a corresponding class,
$Op_{\hbar}(S^{m,k})$, of properly-supported
$\hbar$-pseudodifferential operators acting on $C^{\infty}(M)$.
Moreover, this calculus of operators is independent of the
pariticular choice of partition of unity. There exists a natural
pseudodifferential calculus for such operators with the usual
local symbolic composition formula: Given $a \in S^{m_{1},k_{1}}$
and $b \in S^{m_{2},k_{2}}$, the composition $Op_{\hbar}(a) \circ
Op_{\hbar}(b) = Op_{\hbar}(c) + {\mathcal O}(\hbar^{\infty})$ in
$L^{2}(M)$ where locally,
$$ c(x,\xi;\hbar) \sim \hbar^{-(m_{1}+m_{2})} \sum_{|\alpha|=0}^{\infty}
\frac{{(-i\hbar)}^{|\alpha|} }{\alpha !} (\partial_{\xi}^{\alpha}
a)\cdot (
\partial_{x}^{\alpha} b).$$

\begin{defn} \label{SCQCI}
We say that the operators $Q_{j} \in Op_{\hbar}(S^{m,k}_{cl});
\,\,j=1,...,n$, generate a semiclassical quantum completely
integrable system if
$$ [Q_{i}, Q_{j}]=0; \,\,\, \forall{1 \leq i,j \leq n},$$
\noindent and the respective semiclassical principal symbols
$q_{1},...,q_{n}$ generate a classical integrable system with
$dq_{1} \wedge dq_{2} \wedge \cdot \cdot \cdot \wedge dq_{n} \neq
0$ on a dense open set $\Omega \subset T^{*}M-0$.
\end{defn}

In this article, we are only concerned with Laplacians on compact
Riemannian manifolds. Hence, we only consider the homogeneous case
where the operators $P_{1} = \sqrt{\Delta},P_{2},...,P_{n}$ are
classical pseudodifferential operators of order one. For fixed
$b=(b_{1},b_{2},...,b_{n}) \in {\Bbb R}^{n}$, we define a ladder
of  joint eigenvalues of $P_{1} = \sqrt{\Delta},P_{2},...,P_{n}$
by:
\begin{equation} \label{QCI1}
\{ (\lambda_{1k},...,\lambda_{nk}) \in Spec(P_{1},...,P_{n});
\,\forall j=1,..,n, \,  \lim_{k\rightarrow
\infty}\frac{\lambda_{jk}}{|\lambda_{k}|}  \, = \, b_{j}\},
\end{equation}
\noindent where $|\lambda_{k}|:= \sqrt{ \lambda_{1k}^{2} + ... +
\lambda_{nk}^{2} }$. We denote the corresponding joint
eigenfunctions by  $\phi_{\mu}$. It is helpful to rewrite the
spectral problem in semiclassical notation. The rescaled operators
 $Q_{j}:= \hbar P_{j} \in Op(S^{0,1}_{cl})$ clearly generate a
semiclassical quantum integrable system in the sense of Definition
\ref{SCQCI}. With no loss of generality we may restrict $\hbar$ to
the sequence $\hbar_k = 1/|\lambda_k|$, and then the ladder
eigenvalue problem (\ref{QCI1}) has the form
\begin{equation} \label{QCI3}
Q_{j} \phi_{\mu} =   \mu_{j}(\hbar) \phi_{\mu},
\,\,\,\mbox{where}\,\,\, \mu(\hbar) = b + o(1) \,\,\,\mbox{as}
\,\hbar \rightarrow 0.
\end{equation}
The semi-classical ladder will be denoted by:
\begin{equation} \label{sigma}
 \Sigma_b(\hbar):= \{ \mu(\hbar):=(\mu_{1}(\hbar),...,\mu_{n}(\hbar))
  \in \, \mbox{Spec} (Q_{1},...,Q_{n}); \, | \mu_{j}(\hbar) - b_{j} | \leq C \hbar, \, j=1,...,n \}.
\end{equation}

\subsection{\label{MODELS} Model cases}

Quantum Birkhoff normal forms are microlocal expressions of a
given QCI system in terms of certain model system. Model quantum
completely integrable systems are direct sums of
 the quadratic Hamiltonians:

\begin{itemize}

\item $\hat{I}^{h} := \hbar (D_{y} y + y
D_{y}) $\,\,(hyperbolic \,\, Hamiltonian),

\item $\hat{I}^{e} := \hbar^{2}D_{y}^{2} +
y^{2}$, \,\,(elliptic \,\, Hamiltonian),

\item $ \hat{I}^{ch} := \hbar  \, [ ( y_{1} D_{y_{1}} + y_{2}
D_{y_{2}}) + \sqrt{-1} (y_{1} D_{2}- y_{2}D_{y_{1}}) ] $\,\, (complex \,\, hyperbolic \,\, Hamiltonian),

\item $ \hat{I} := \hbar \,  D_{\theta } $, \,\, (regular \,\, Hamiltonian).

\end{itemize}

\vspace{2mm}

\noindent The corresponding  model eigenfunctions are:

\begin{itemize}

\item  $u_{h}(y; \lambda, \hbar) =  |\log \hbar|^{-1/2} \, [
c_{+}(\hbar)  Y(y)\,|y|^{-1/2 + i\lambda(\hbar)/\hbar} +
c_{-}(\hbar)
 Y(-y) \, |y|^{-1/2 + i\lambda(\hbar)/\hbar} ];\, |c_{-}(\hbar)|^{2} + |c_{+}(\hbar)|^{2} =1; \, \, \lambda(\hbar) \in {\Bbb R}$.

\item $u_{ch}(r, \theta ; t_{1}, t_{2}, \hbar) =  |\log \hbar|^{-1/2}
r^{(-1 + it_{1}(\hbar) )/\hbar} \,e^{it_{2}(\hbar) \theta } ; \,\,
t_{1}(\hbar), t_{2}(\hbar) \in {\Bbb R}.$

\item $u_{e}(y;n, \hbar) = \hbar^{-1/4} \exp ( -
y^{2}/\hbar ) \,\,\Phi_{n}(\hbar^{-1/2}y) ; \,\, n \in {\Bbb N}.$

\item $u_{reg}( \theta ;m,  \hbar) = e^{i m \theta}; \,\, m \in {\Bbb Z}.$

\end{itemize}

\noindent Here, $Y(x)$ denotes the Heaviside function,
$\Phi_{n}(y)$ the $n$-th Hermite polynomial  and $(r,\theta )$
polar variables in the $(y_{1}, y_{2})$ complex hyperbolic plane.
We restrict to those model eigenfunctions with  $\lambda(\hbar),
t_{1}(\hbar), t_{2}(\hbar) = {\mathcal O}(\hbar)$ as $\hbar
\rightarrow 0$, since these are the ones which localize properly.
The important part of a model eigenfunctions is its
microlocalization to a neighborhood of $x = \xi = 0$, so we put:
$$ \psi(x;\hbar):= Op_{\hbar}( \, \chi(x) \chi(y) \chi(\xi) \, )  \cdot  u(y;\hbar),$$
where $\epsilon >0$ and  $\chi \in C^{\infty}_{0}([-\epsilon,
\epsilon])$. In the hyperbolic, complex hyperbolic, elliptic  and
regular cases, we write  $\psi_{h}(y;\hbar),  \psi_{ch}(y ;\hbar),
\psi_{e}(y ;\hbar)$ and $ \psi_{reg}(y ;\hbar)$ respectively. A
straightforward computation \cite{T2} shows that when
$t_{1}(\hbar), t_{2}(\hbar), n \hbar , m \hbar = {\mathcal
O}(\hbar)$  the model quasimodes are $L^{2}$-normalized; that is
\begin{equation} \label{OP} \| Op_{\hbar}(, \chi(x) \chi(y) \chi(\xi) \, ) \, u(y;\hbar)
\|_{L^{2}} \sim 1 \end{equation}  as $\hbar \rightarrow 0$. Note
that, although the model eigenfunctions above are not in general
smooth functions, the microlocalizations are $C^{\infty}$ and
supported near the origin.

\subsection{Singular Birkhoff normal form}

In this section we introduce Birkhoff normal forms for a quantum
completely integrable system near a singular orbit. The main result
 is the following  quantum analogue (see also \cite{VN}) of the classical Eliasson normal form in (\ref{E5}):

\begin{lem} \label{QBNF}
 Let $c \in {\Bbb R}^{n}$ be a singular value of the moment map ${\mathcal P}$ and ${\Bbb R}^{n} \cdot v_{0}$ be a rank-$k$ Eliasson non-degenerate orbit of the joint flow. Then, there exists a microlocally
elliptic $\hbar$-Fourier integral operator, $F_{\hbar}$, and a microlocally
invertible $n\times n$ matrix of pseudodifferential operators,
${\mathcal M}_{ij},$ with $[ M_{ij}, Q_{k}]=0; k=1,...,n$ and satisfying:
\begin{equation} \label{E6}
F_{\hbar}^{-1} \, ( Q_{1} - c_{1},...,Q_{n} - c_{n} ) \, F_{\hbar}  =_{ \Omega}
 {\mathcal M} \,  \cdot  ( \, \hat{I}^{h} - \epsilon^{h}(\hbar), \hat{I}^{ch} - \epsilon^{ch}(\hbar), \hat{I}^{e} - \epsilon^{e}(\hbar), \, \hbar D_{\theta} - \epsilon(\hbar) \, )  + {\mathcal O}(\hbar^{\infty}).
\end{equation}

\noindent Here, $\hat{I}^{h}_{j} = \hbar (D_{y_{j}} y_{j} + y_{j}
D_{y_{j}}), \,\hat{I}^{e}_{j} = \hbar^{2}D_{y_{j}}^{2} +
y_{j}^{2}, \,\hat{I}^{ch} = \hbar [ ( y_{j} D_{y_{j}} + y_{j+1}
D_{y_{j+1}}) + \sqrt{-1} (y_{j} D_{j+1}- y_{j+1}D_{y_{j}}) ] $ and $\hat{I}_{j} = \hbar D_{\theta_{j}}.$
\end{lem}

\noindent{\em Proof:} The proof is essentially the same as in
([VN] Theorem 3.6). The only complication here is that since
${\Bbb R}^{n} \cdot(v_{0})$ is a rank $k<n$ torus and not a point,
$\hbar D_{\theta_{1}},...,\hbar D_{\theta_{k}}$ must be added to
the space of model operators. The proof can be reduced to that in
[VN] by making Fourier series decompositions in the
$(\theta_{1},...,\theta_{k})$ variables (see, for instance \cite{T2}
Theorem 3). \qed

 Consider the microlocal (quasi-)eigenvalue problem in $\Omega:$ \begin{equation}
\label{QEP}  Q_{j} \psi_{\nu} =_{\Omega} \nu_{j}(\hbar)
\psi_{\nu}; j=1,...,n. \end{equation} We say that
$\nu_k(\hbar)$ is a quasi-classical eigenvalue if there exists a
non-trivial solution of (\ref{QEP}). The set of quasi-classical
eigenvalues around $c$ is thus:
\begin{equation} \label{sigma}
 Q\Sigma_c(\hbar):= \{ \nu (\hbar): (\ref{QEP}) \; \mbox{holds, with }\;  \,
  | \nu (\hbar) - c | \leq C \hbar, \, j=1,...,n \}.
\end{equation}

The solution space of (\ref{QEP}) can be characterized  uniquely
(up to a ${\Bbb C}(\hbar)$-multiple) in terms of the model
quasimodes $\psi_{e}, \psi_{h}, \psi_{ch}$ and $\psi_{reg}$. In
the following, we use the abbreviation   $( u_{e}(y;n, \hbar)
\cdot u_{h}(y;\lambda_{k}(\hbar), \hbar) \cdot
u_{ch}(y;t_{1,k}(\hbar),t_{2,k}(\hbar), \hbar) \prod_{j=1}^{k}
e^{im_{j}\theta_{j}} \, )$   for the expression $$\Pi_{j = 1}^H
\psi_{h j} (y_j;\hbar) \otimes \Pi_{j = H + 1}^{H + L} \psi_{ch j}
(y_j;\hbar) \otimes \Pi_{j = H + L+ 1}^{H + L + E} \psi_{e j}
(y_j;\hbar) \otimes \Pi_{j = L + H + E + 1}^{n} \psi_{r j}
(y_j;\hbar)$$ in which the $(m, n, \lambda, t_1, t_2)$ parameters
are put in.

\begin{prop} \label{FORMULA} For any admissible solution $\psi_{\nu (\hbar)}$ of (\ref{QEP}), there
 exist
$t_{1}(\hbar),t_{2}(\hbar), \lambda (\hbar) \in \R$ and $
n, m \in {\bf N}$ and an $\hbar$-dependent
constant $c(\hbar)$ such that
$$ \psi_{\nu} =_{\Omega} c(\hbar) \; F \;\;( u_{e}(y;n, \hbar) \cdot
u_{h}(y;\lambda (\hbar), \hbar) \cdot
u_{ch}(y;t_{1}(\hbar),t_{2}(\hbar), \hbar) \prod_{j=1}^{k}
e^{im_{j}\theta_{j}} \, ) ,$$ where $F_{\hbar}$ is the microlocally
unitary $\hbar$-Fourier integral operator in Lemma (\ref{QBNF}).
\end{prop}

\begin{proof}
The proposition follows from the uniqueness of microlocal solutions
 of the model eigenfunction equations (see \cite{CP} and \cite{VN})
 up to ${\Bbb C}(\hbar)$-multiples and  Lemma (\ref{QBNF}). 
\end{proof}

\section{Blow-up of eigenfunctions attached to singular leaves of the
Lagrangian fibration}

The purpose of this section is to  prove Theorem (\ref{LP}) and
Lemma (\ref{BLOWUP}). The key estimate is a small scale $L^2$ mass
estimate near a singular leaf $\Lambda$. Roughly speaking, it says
that the $L^{2}$-mass of the eigenfunctions $\phi_{\mu}$
with $\mu(\hbar) \in \Sigma_c(\hbar)$  in a tube of radius
$\hbar^{\delta}$ around $\pi(\Lambda)$ satisfies
\begin{equation} \label{dom.5}
  ( \,  Op_{\hbar}( \chi_{1}^{\delta} (x;\hbar)) \phi_{\mu}, \phi_{\mu}  \, )
  \geq C > 0,
\end{equation}
\noindent where $C$ is a constant independent of $\hbar$ and where
 $\chi_{1}^{\delta}(x;\hbar) := \chi_{1}(\hbar^{-\delta}x)$ with
$\chi_{1}(x)$  a cutoff supported near $\pi (\Lambda)$.

 In the proof, we will need to use additional
pseudo-differential operators belonging to  a more refined
 semi-classical calculus, containing cutoffs such as $\chi_{1}(\hbar^{-\delta}x)$, which
involve smaller length scales.

\subsection{Eigenfunction mass near non-degenerate, singular orbits}

We first give a lower bound for the microlocal mass on `large'
length scales of joint eigenfunctions with joint eigenvalues
$\mu(\hbar) \in \Sigma_c(\hbar)$ near $\Gamma_{sing}(c)$. Let
$\chi_{reg} \in C^{\infty}_{0}(T^{*}M)$ and $\chi_{sing} \in
C^{\infty}_{0}(T^{*}M)$ be cutoff functions supported near
$\Gamma_{reg}(c)$ and $\Gamma_{sing}(c)$ respectively.  We claim
 that for any $f \in S(R)$ with $\check{f} \in
C^{\infty}_{0}({\Bbb R})$ with sufficiently small support and
$\hbar \in (0, \hbar_{0}],$ we have
\begin{equation} \label{mass3}
\sum_{j=1}^{\infty} \langle Op_{\hbar}(\chi_{sing}) \phi_{\mu},
\phi_{\mu} \rangle \, f ( \hbar^{-1}( \mu_j(\hbar) - c) ) \, \geq
\, \frac{1}{C_{0}} >0.
\end{equation}
 Here, the sum runs over the entire joint spectrum.  The
estimate (\ref{mass3}) follows readily from a
local (singular) Weyl law near a singular rank-$k$ orbit ${\Bbb R}^{n} \cdot v_{0}.$ The local Weyl law inturn follows from the Eliasson normal form and standard wave trace computations (\cite{BPU} Theorems 1.1-1.3). When the singular orbit has only hyperbolic summands, the RHS in (\ref{mass3}) can be improved to  $|\log \hbar |^{k}$ (see \cite{BPU} ).

Consequently, for each $\hbar$, there exists at least one  joint
eigenfunction $\phi_{\mu}$ of the $Q_{j}$'s with joint
eigenvalue $\mu (\hbar) \in \Sigma_c (\hbar)$ satisfying:

\begin{equation} \label{mass4}
\langle Op_{\hbar}(\chi_{sing}) \phi_{\mu}, \phi_{\mu} \rangle \geq \frac{1}{C_{0}} >0,
\end{equation}

provided we choose the constant $C>0$ in (\ref{sigma}) large
enough. By the local Weyl law  in the regular case, again for $C>0$ large enough in
(\ref{sigma}),  there also exist joint eigenfunctions $\phi_{\mu}$ with $\mu = \mu (\hbar) \in \Sigma (\hbar)$ satisfying

\begin{equation} \label{mass4.1}
\langle Op_{\hbar}(\chi_{reg}) \phi_{\mu}, \phi_{\mu} \rangle \geq  \frac{1}{C_{0}}.
\end{equation}
 Here, $C_{0}>0$ denotes possibly different constants from line to line.

\begin{defn}
Let $V_{sing}(\hbar) \, (resp. V_{reg}(\hbar) )$ denote the set of $L^{2}$-normalized
 joint eigenfunctions of $Q_{1},...,Q_{n}$ with $  \mu(\hbar)\in  \Sigma_c(\hbar)$  and satisfying (\ref{mass4}) \, (resp. \ref{mass4.1}).
\end{defn}

 An easy argument involving invariant measures (see for example \cite{TZ1})
 shows that the connected component $\Gamma_{sing}(c)$ must contain a rank-$k <n$ compact
  orbit ${\Bbb R}^{n}\cdot (v)$. Again, to simplify the writing we can, without loss of
  generality that there is precisely one such compact orbit in $\Gamma_{sing}(c)$. As the
  next Lemma shows,  the mass much concentrate on at least one of
  these singular orbits $\Lambda$.

\medskip

\begin{lem} \label{LAMBDA}
Let  $\Lambda:= {\Bbb R}^{n} \cdot(v_{0})$ be an Eliasson
non-degenerate orbit, and let $\chi_{\Lambda}(x,\xi) \in
C^{\infty}_{0}(\Omega;[0,1])$ be a cutoff function supported in an
invariant  neighbourhood $\Omega$ of $\Lambda$  and identically
equal to one on a smaller neighbourhood $\tilde{\Omega}$ of
$\Lambda$. Then there exists a constant
 $C_{0}>0$ such that for $\hbar \in
(0,\hbar_{0}]$ and $\phi_{\mu} \in V(\hbar)$,
 $$( Op_{\hbar}(\chi_{\Lambda}) \phi_{\mu}, \phi_{\mu} ) \geq \frac{1}{C_{0}}.$$
\end{lem}
\noindent{\em Proof}:  With no essential loss of generality, we
will assume for notational simplicity that
$${\mathcal P}^{-1}(c) = \, \Gamma_{sing}(c) \cup \, \Gamma_{reg}(c)$$
\noindent and that $\Gamma_{sing}(c)$ contains a single compact
rank-$k$ orbit, $ \Lambda:= {\Bbb R}^{n}(v_{0})$.  The argument
easily extends to the case of multiple orbits and components.

 Let
 $  \chi_{sing} \in C^{\infty}_{0}(T^{*}M)$ be the
 cutoff functions defined above (cf. \ref{mass3} - \ref{mass4}).
Since $\chi_{\Lambda} \chi_{sing} = \chi_{\Lambda}$, hence
$Op_{\hbar}(\chi_{\Lambda} ) \circ Op_{\hbar}(\chi_{sing}) =
Op_{\hbar} (\chi_{\Lambda}) + O(\hbar),$ it follows from
(\ref{mass4}) that
\begin{eqnarray} \label{mass6}
1/C \leq   (Op_{\hbar}(\chi_{\Lambda}) \phi_{\mu}, \phi_{\mu}) + (
(1- Op_{\hbar}(\chi_{\Lambda}) ) \circ Op_{\hbar}(\chi_{sing})
\phi_{\mu}, \phi_{\mu}).
\end{eqnarray}

 \noindent Since non-compact orbits do not support probability measures which are invariant under
 the joint flow $\Phi_{t}$ (see \cite{TZ1}), by the compactness of ${\mathcal P}^{-1}(c)$ it follows that $ \Lambda = {\Bbb R}^{n} \cdot(v_{0})$,
the compact  rank $ k < n$ orbit of $v_{0} \in {\mathcal
P}^{-1}(c)$, is the forward limit set for the joint flow on
${\mathcal P}^{-1}(c)$. The  semiclassical Egorov Theorem then
gives:
$$ ( (1- Op_{\hbar}(\chi_{\Lambda}) ) \circ Op_{\hbar}(\chi_{sing}) \phi_{\mu}, \phi_{\mu}) =
 ( Op_{\hbar}(\Phi_{t}^{*}(1- \chi_{\Lambda}) ) \circ Op_{\hbar}(\chi_{sing}) \phi_{\mu},
\phi_{\mu}) + {\mathcal O}(\hbar).$$

\noindent Moreover, since $ \Lambda = {\Bbb R}^{n} \cdot(v_{0})$
is an $\omega$-limit set for the joint flow $\Phi_{t}$ on
${\mathcal P}^{-1}(c)$, there exists $t_{0} \in \Bbb R^{n}$ such
that: \vspace{3mm}
\begin{itemize}

\item $ supp ( \Phi_{t_{0}}^{*}(1- \chi_{\Lambda}) ) \cap supp ( \chi_{sing}) \subset  \Omega$, \,\,\mbox{where the  QBNF is valid in} $\Omega$.

\vspace{2mm}

\item $ supp ( \Phi_{t_{0}}^{*}(1- \chi_{\Lambda}) ) \cap {\Bbb R}^{n} \cdot(v_{0}) = \emptyset.$

\end{itemize}
\vspace{3mm}
  In ([T2] Lemma 6), it is shown that microlocal eigenfunctions $\phi_{\mu}$ given by the QBNF in Lemma (\ref{QBNF}) near
  the orbit $\Lambda$ satisfy  \begin{equation} \label{OFFDIAGONAL}  (
Op_{\hbar}(\Phi_{t_{0}}^{*}(1- \chi_{\Lambda}) ) \circ
Op_{\hbar}(\chi_{sing}) \phi_{\mu}, \phi_{\mu}) = {\mathcal
O}(|\log \hbar|^{-1}). \end{equation}  The lemma then follows from
(\ref{mass6}),  (\ref{OFFDIAGONAL}).\qed

\medskip

\subsection{Localization on singular orbits}

 We claim that the joint eigenfunctions $ \phi_{\mu} \in V(\hbar)$ must blow
  up along $\pi(\Lambda)$. The first way of quantifying this blowup
   involves computing the asymptotics for the expected values $( Op_{\hbar} (q) \phi_{\mu}, \phi_{\mu})$
   where $q \in C^{\infty}_{0}(T^{*}M)$.
\begin{lem} \label{C} Let  $ \phi_{\mu} \in V(\hbar)$. Then:

\begin{equation} \label{blowup1}
(Op_{\hbar}(q) \phi_{\mu}, \phi_{\mu}) = |c(\hbar)|^2 \int_{{\Bbb
R}^{n} \cdot(v_{0})} q \, d\mu + {\mathcal O}(|\log
\hbar|^{-1/2}),
\end{equation}
\noindent again with $|c(\hbar)|  \geq \frac{1}{C_{0}}.$
 \end{lem}

\begin{proof} Since $\phi_{\mu}$ solves the equation (\ref{QEP})
exactly  (and afortiori microlocally on $\Omega$), we may express
it by Proposition (\ref{FORMULA})  in the form:
\begin{equation}\label{FORM2}
\phi_{\mu} =_{\Omega} c(\hbar)  \; F \; u_{\mu},
\end{equation}
for some constant $c(\hbar)$. Here, \begin{equation} \label{U}
u_{\mu} = ( u_{e}(y;n, \hbar) \cdot u_{h}(y;\lambda (\hbar),
\hbar) \cdot u_{ch}(y;t_{1}(\hbar),t_{2}(\hbar), \hbar)
\prod_{j=1}^{k} e^{im_{j}\theta_{j}} \, )
\end{equation} Here, by applying the operators on both sides of the QBNF in Lemma (\ref{QBNF}) to the model eigendistributions $u_{\mu}$ and using the uniqueness  result in Lemma (\ref{FORMULA}), it follows that for some $n \times n$ matrix $M$ with $M(0) \in GL_{n},$
$$M (m\hbar , n \hbar , \lambda(\hbar), t_{1}(\hbar), t_{2}(\hbar) ) \, \cdot  \, ( m\hbar , n \hbar , \lambda(\hbar), t_{1}(\hbar), t_{2}(\hbar) ) = \mu(\hbar).$$
By the inverse function theorem, the $(m\hbar , n \hbar , \lambda(\hbar), t_{1}(\hbar), t_{2}(\hbar))$ are uniquely determined (modulo ${\mathcal O}(\hbar^\infty) $ ) by the joint eigenvalues $\mu(\hbar)$ and  moreover, when $\mu(\hbar) \in \Sigma(\hbar)$ it follows that $ m\hbar , n \hbar , \lambda(\hbar), t_{1}(\hbar), t_{2}(\hbar) = {\mathcal O}(\hbar).$ 
By Lemma
(\ref{LAMBDA}), by (\ref{OP})  and by (\ref{FORM2}), it follows
that for $\hbar \in (0,\hbar_{0}],$
\begin{equation} \label{CH} \begin{array}{l}
\frac{1}{C_{0}} \leq ( Op_{\hbar}(\chi_{\Lambda}) \phi_{\mu},
\phi_{\mu} ) = c(\hbar) ( F^* Op_{\hbar}(\chi_{\Lambda}) F
u_{\mu}, u_{\mu} ) \leq \;\; |c(\hbar)|.
\end{array}
\end{equation}

Granted this lower bound on $|c(\hbar)|$, the Lemma reduces to
estimating matrix elements of model eigenfunctions. We now
evaluate the matrix elements case by case. The most interesting
case is where the orbit $\Lambda$ is strictly real or complex
hyperbolic. Then, given any $q \in C^{\infty}_{0}(T^{*}M)$, it
follows from (\ref{OFFDIAGONAL}) that
\begin{equation}
(Op_{\hbar}(q) \phi_{\mu}, \phi_{\mu}) = (Op_{\hbar}(q) \circ
Op_{\hbar}(\chi_{\Lambda}) \phi_{\mu}, Op_{\hbar}(\chi_{\Lambda})
\phi_{\mu}) + {\mathcal O}(|\log \hbar|^{-1}).
\end{equation}
\noindent We now use (\ref{FORM2}) to conjugate to the model
setting. The function $q$ goes to $q \circ \chi$ where $\chi$ is
the canonical transformation underlying $F$. The model $\R^n$-
action locally reduces to a compact torus $T^k$-action, so we can
average the function $q \circ \chi$ over the action to obtain a
smooth invariant function. We then Taylor expand this averaged
function in the directions $(y, \eta)$  transverse to the action.
We obtain:
\begin{equation}
(Op_{\hbar}(q) \circ Op_{\hbar}(\chi_{\mu}) \phi_{\mu},
Op_{\hbar}(\chi_{\Lambda}) \phi_{\mu}) = |c (\hbar)|^2
\int_{{\Bbb R}^{n} \cdot(v_{0})} q \, d\mu +
(Op_{\hbar}(r_{h})u_{\mu}, u_{\mu}) + (Op_{\hbar}(r_{ch})u_{\mu},
u_{\mu}) + {\mathcal O}(\hbar).
\end{equation}
\noindent \noindent where, $r_{h}, r_{ch} \in
C^{\infty}_{0}(\Omega)$ are the Taylor remainders  with $r_{h},
r_{ch}= {\mathcal O}(|y| + |\eta|)$.  A direct computation for the
model distributions, $u_{\mu}$ (see [T2] Lemma 5 and
Proposition 3) shows that:
\begin{equation}
(Op_{\hbar}(r_{h})u_{\mu}, u_{\mu}) = {\mathcal O}(|\log
\hbar|^{-1/2}), \,\,\,(Op_{\hbar}(r_{ch})u_{\mu}, u_{\mu}) =
{\mathcal O}(|\log \hbar|^{-1/2}).
\end{equation}

The remaining cases are where elliptic (i.e. Hermite factors).
Each such factor satisfies
\begin{equation} \label{blowup2}
(Op(r_{e})u_{\mu}, u_{\mu}) = {\mathcal O}(\hbar) \end{equation} so is
better than what is claimed.
 \end{proof}

\subsection{Lower bounds on $L^p$ norms: Proof of Lemma 0.3}

Our objective in this section is to refine  the argument in Lemma
3.3 to actually produce pointwise  lower bounds for eigenfunctions
attached to Eliasson non-degenerate leaves of the Lagrangian
fibration.  To do this, we study matrix elements on much smaller
length scales. That requires the use of small scale
pseudodifferential operators. We pause to define these objects.

\subsubsection{Small scale semiclassical pseudo-differential
calculus}

The more refined symbols are defined as follows: Given an open set
$U \in {\Bbb R}^{n}$ and $0 \leq \delta < \frac{1}{2}$, we say
that $a(x,\xi;\hbar) \in S^{m}_{\delta}(U\times {\Bbb R}^{n})$if
\begin{equation}
| \partial_{x}^{\alpha} \partial_{\xi}^{\beta} a(x,\xi;\hbar) |
\leq C_{\alpha \beta} \hbar^{-\delta(|\alpha |+|\beta |)}.
\end{equation}
Model symbols include cutoffs of the form $\chi(h^{-\delta} x,
h^{-\delta} \xi)$ with $\chi \in C_0^{\infty}(\R^{2n})$.  There is
a pseudodifferential calculus $Op_{\hbar} S^{m}_{\delta}(U\times
{\Bbb R}^{n})$ associated with such symbols with the usual
symbolic composition formula and Calderon-Vaillancourt
$L^{2}$-boundedness theorem [Sj]. Composition with operators  in
our original class $Op_{\hbar} S^{m,0}(U \times {\Bbb R}^{n})$
preserves $Op_{\hbar} S^{m}_{\delta}(U\times {\Bbb R}^{n})$.

\subsubsection{Outline of proof}

Let us now  outline the proof. If $\chi_{2}^{\delta}(x,\xi;\hbar)
\in C^{\infty}_{0}(T^{*}M)$ is a
 second cutoff supported in a radius $\hbar^{\delta}$ tube around $\Lambda$ then clearly
\begin{equation} \label{dom1}
\chi_{1}^{\delta} \geq \chi_{2}^{\delta}
\end{equation}
\noindent where we view both functions in (\ref{dom1}) as being defined on $T^{*}M$ and so (\ref{dom1})
 is just a pointwise estimate. Modulo small errors (see (\ref{WG})),   inequality (\ref{dom1})
 implies the corresponding operator bound  for the matrix elements:
\begin{equation} \label{dom2}
 ( \, Op_{\hbar}(\chi^{\delta}_{1}) \phi_{\mu}, \phi_{\mu} \, ) \gg ( \, Op_{\hbar}(\chi^{\delta}_{2}) \phi_{\mu}, \phi_{\mu} \, ).
\end{equation}

 Now, take $\phi_{\mu} \in V(\hbar)$. From the expressions (\ref{FORM2}) and
 (\ref{CH}), it follows that
 the matrix elements on the RHS of (\ref{dom2})
   can be estimated from below by simply computing the masses
   of the model distributions $u(y,\theta;\hbar) = \prod_{j=1}^{k} e^{im_{j}\theta_{j}}
   u_{e}(y) \cdot u_{h}(y) \cdot u_{ch}(y)$. Our purpose now is to compute these masses  on {\em shrinking}
    neighbourhoods of diameter $\hbar^{\delta}$ centered around a singular rank $\ell <n$ orbit, $\Lambda$. We
    show that in such neighbourhoods of diameter $\hbar^{\delta}$,  these model distributions have  finite
    mass bounded from below by a positive constant independent of $\hbar \in (0,\hbar_{0}]$ provided we
    choose  $\delta < 1/2$. It  follows that,  for $\delta <1/2$,
\begin{equation}\label{dom3}
(Op_{\hbar}(\chi_{2}^{\delta}) \phi_{\mu}, \phi_{\mu}) \geq
\frac{1}{C}
>0.
\end{equation}

Then, combined with (\ref{dom2}), the lower bound (\ref{dom3}) implies the $L^{\infty}$ lower bound in Lemma 0.3.
 We now turn to the details of the proof:

\vspace{2mm}

\begin{proof}

By choosing a sufficiently small tubular neighbourhood, $\tilde{M}
\supset M$ inside $T^{*}M$, we can extend the Riemannian metric on
$M$  to a Riemannian metric on the entire tube, $\tilde{M}$.
Moreover, by possibly rescaling $\hbar$ by an appropriate
constant, we can assume that ${\mathcal P}^{-1}(c) \subset
\tilde{M}$. We will make this assumption without further comment.
We will continue to denote by $A$ and $A_{\epsilon}$ the tubular
neighbourhoods defined above. In addition, we define
\begin{equation}
A_{\hbar^{\delta}} = \exp \{ (x,v) \in E(\delta); |v| \leq
\hbar^{\delta} \}.
\end{equation}

\noindent Let $\chi_{1}^{\delta} (x;\hbar ) \in
C^{\infty}_{0}(A_{\hbar^{\delta}};[0,1])$ be a cutoff function
which is identically equal to $1$ in $A_{\hbar^{\delta}/2}$.
Clearly,
\begin{equation}
\chi_{1}^{\delta}(x;\hbar) \in S^{0}_{\delta}.
\end{equation}

\noindent Now, choose another cutoff function,
$\chi_{2}^{\delta}(x,\xi;\hbar)) \in C^{\infty}_{0}(T^{*}M)$,
which is supported in neighbourhood of radius $\hbar^{\delta}$
around the torus $\Lambda = {\Bbb R}^{n} \cdot(v_{0})$, with the
property that $\chi_{1}^{\delta} \geq \chi_{2}^{\delta}$
pointwise. By the Garding inequality, there exists a constant
$C_{1}>0$ such that:
\begin{equation} \label{WG}
(Op_{\hbar}(\chi_{1}^{\delta}) \phi_{\mu},\phi_{\mu)}) \geq
(Op_{\hbar}(\chi_{2}^{\delta})\phi_{\mu}, \phi_{\mu}) - C_{1}
\hbar^{1-2\delta}.
\end{equation}

We now conjugate the right side to the model by the $\hbar$
Fourier integral operator $F$ of Lemma (\ref{QBNF}). With no loss
of generality, we may assume the cutoff function
$\chi_{1}^{\delta}(y,\eta,I;\hbar)$ to be of product type:
$$ \chi_{1}^{\delta}(y,\eta,I;\hbar) = \prod_{j=1}^{L+N}
\chi(\hbar^{-\delta} y_{j}) \, \chi(\hbar^{-\delta} \eta_{j})
\cdot \prod_{j=L+N+1}^{L+M+N+1} \chi(\hbar^{-\delta}\rho_{j})
\chi(\hbar^{-\delta} \alpha_{j}) \cdot \prod_{j=n-l+1}^{n}
\chi(\hbar^{-\delta} I_{n+1-j}).$$ \noindent Here $(r_{j},
\alpha_{j})$ denote radial variables in the $j$-th complex
hyperbolic summand.  Since $F$ is a microlocally elliptic
$\hbar$-Fourier integral operator associated to a canonical
transformation $\kappa$, it follows by Egorov's theorem
\begin{equation} \label{bound}
(Op_{\hbar}(\chi_{2}^{\delta}) \phi_{\mu}, \phi_{\mu}) =
|c(\hbar)|^2  (Op_{\hbar}(\chi_{2}^{\delta} \circ \kappa) u_{\mu},
u_{\mu}) - C_{3} \hbar^{1-2\delta}
\end{equation}
\noindent where $c(\hbar) u_{\mu} (y,\theta;\hbar)$ is the
microlocal normal form (\ref{U}) for the eigenfunction
$\phi_{\mu}$. To simplify the notation a little, we will write
$\chi^{\delta}(x;\hbar) := \chi(\hbar^{-\delta} x)$ below. Now,
$(\tilde{\chi}^{\delta} u, u)$ consists of products of four types
of terms. The first three are:
$$M_{e} = \int_{-\infty}^{\infty} \chi^{\delta}(\eta;\hbar) \,|
\widehat{ {\chi}^{\delta} u_{e}} (\eta;\hbar ) |^{2} d\eta,$$
$$M_{h} = \int_{-\infty}^{\infty} \chi^{\delta}  (\eta;\hbar) \,|
\widehat{ \chi^{\delta} u_{h} }(\eta;\hbar) |^{2} d\eta,$$
\noindent and finally,
$$M_{ch} = \int_{-\infty}^{\infty} \int_{-\infty}^{\infty}\chi^{\delta}
(\eta_{1},\eta_{2};\hbar ) \,| \widehat{ \chi^{\delta} u_{ch}
}(\eta_{1},\eta_{2};\hbar) |^{2}  \,d\eta_{1} d\eta_{2} .$$
\noindent To estimate $M_{e}$, we note that, since
$\phi_{\mu} \in V(\hbar),$ and
$$ {\mathcal F}(e^{-|y|^{2}/\hbar} \Phi_{n}(y \hbar^{-1/2}))(\eta) =
e^{-|\eta|^{2}/\hbar} \Phi_{n}(\eta \hbar^{-1/2}),$$ \noindent it
follows that,
$$M_{e} \gg \int_{-\infty}^{\infty} e^{-2|\eta|^{2}/\hbar} |\Phi_{n}(\hbar^{-1/2}
\eta)|^{2} d\eta + {\mathcal O}(\hbar^{\infty})$$ \noindent and so
for $\hbar \in (0, \hbar_{0}], \, M_{e}(\hbar) \gg 1 $.

To estimate $M_{h}$, we write:
\begin{equation} \label{modelmass1}
M_{h} = \frac{  |c_{+}(\hbar)|^{2} + |c_{+}(\hbar)|^{2}} { \log
\hbar} \, \left( \int_{0}^{\infty} \chi(\hbar \xi/\hbar^{\delta})
\left| \int_{0}^{\infty} e^{-ix} x^{-1/2 + i \lambda/\hbar}
\chi(x/\hbar^{\delta}\xi) dx \right|^{2} \frac{d\xi}{\xi} \right).
\end{equation}
\noindent Since, $ \phi_{\mu} \in V(\hbar),$ it follows that $|c_{+}(\hbar)|^{2} + |c_{+}(\hbar)|^{2} \gg 1$ and  the integral in  (\ref{modelmass1}) is bounded  from below by
\begin{equation} \label{modelmass2}
\frac{1}{C_{0}} \, (\log \hbar)^{-1} \, \int_{0}^{\hbar^{\delta
-1} } \frac{d\xi}{\xi} \left| \int_{0}^{\hbar^{\delta} \xi}
e^{-ix} x^{-1/2 + i \lambda/\hbar} dx \right|^{2} + {\mathcal
O}(|\log \hbar|^{-1}).
\end{equation}
\noindent To estimate this last integral, assume first that $\xi
\in [0,\hbar^{-\delta}]$.  Then, by an integration by parts,
$$\int_{0}^{\hbar^{\delta} \xi}  e^{-ix} x^{-1/2 + i \lambda/\hbar} dx =
{\mathcal O}(|\hbar^{\delta}\xi|^{1/2}) + {\mathcal
O}(|\hbar^{\delta}\xi|^{3/2}),$$ \noindent and so,
\begin{equation} \label{modelmass3}
M_{h} \gg |\log \hbar|^{-1} \int_{\hbar^{-\delta}}^{\hbar^{\delta
-1} } \frac{d\xi}{\xi} \left| \int_{0}^{\hbar^{\delta} \xi}
e^{-ix} x^{-1/2 + i \lambda/\hbar} dx \right|^{2} + {\mathcal
O}(|\log \hbar|^{-1}).
\end{equation}
\noindent From (\ref{modelmass3}), it follows that:
\begin{equation} \label{modelmass4}
M_{h} \gg \,  |\Gamma (1/2 + i \lambda/\hbar)|^{2}\, (1-2\delta) +
{\mathcal O}(|\log \hbar|^{-1}).
\end{equation}
\noindent Since, $\delta = 1/2-\epsilon$, there exists a constant
$C(\epsilon)$ such that,  $M_{h} \geq C(\epsilon) >0$ uniformly
for $\hbar \in (0,\hbar_{0}(\epsilon)]$.

Finally, we are left with the integral $M_{ch}$ corresponding to a
loxodromic subspace. Since $|{\mathcal J}_{k}(\rho)| \leq 1$ for
all $k \in {\Bbb Z}$ and $\rho \in {\Bbb R}$, it follows that:
\begin{equation} \label{modelmass5}
M_{ch} \gg |\log \hbar|^{-1} \int_{\hbar^{-\delta}}^{\hbar^{\delta
- 1}} \left| \int_{0}^{\hbar^{\delta}\alpha} {\mathcal
J}_{k}(\rho) \rho^{it/\hbar} d\rho \right|^{2}
\frac{d\alpha}{\alpha} + {\mathcal O}(|\log \hbar|^{-1}).
\end{equation}
\noindent Here, ${\mathcal J}_{k}(\rho)$  denotes the $k$-th
integral Bessel function of the first kind [AS]. For $\alpha \geq
\hbar^{-\delta}$,

\begin{equation} \label{modelmass6}
| \int_{0}^{\hbar^{\delta}\alpha} {\mathcal J}_{k}(\rho)
\rho^{it/\hbar} d\rho | = \left|  \frac{2^{it/\hbar} \Gamma
\frac{k+1+it/\hbar}{2})}{\Gamma( \frac{k+1 -it/\hbar}{2}) }
\right| + {\mathcal O}( |\hbar^{\delta}\alpha|^{-1/2}) = 1+
{\mathcal O}( |\hbar^{\delta}\alpha|^{-1/2}), \vspace{-5mm}
\end{equation}

 and so,
\begin{equation} \label{modelmass7}
M_{ch} \gg |\log \hbar|^{-1} \int_{ {\hbar}^{-\delta} }^{
\hbar^{\delta - 1}} \frac{d\alpha}{\alpha} + {\mathcal O}(|\log
\hbar|^{-1}) = 1-2\delta  + {\mathcal O}(|\log \hbar|^{-1}).
\end{equation}
Consequently, given $\delta = 1/2 -\epsilon$ it again follows that
$M_{ch} \geq C(\epsilon) >0$ uniformly for $\hbar \in
(0,\hbar_{0}(\epsilon)]$.

The final step involves estimating $( Op_{\hbar}(\chi^{\delta}
(I)) e^{im\theta}, e^{im\theta} )$. An integration by parts in the
$I_{1},...,I_{\ell}$ variables shows that:
\begin{equation}
( Op_{\hbar} ( \chi^{\delta} (I) ) e^{im\theta}, e^{im\theta} ) =
1 + {\mathcal O}(\hbar^{1-\delta}).
\end{equation}
As a consequence of the estimates above for $M_{h}, M_{ch}, M_{e}$ and the bound in (\ref{bound}) it follows that for {\em any} $\epsilon >0$ and  $\delta =
1/2-\epsilon$, there exists a constant $C(\epsilon)>0$ such that
for all $\phi_{\mu} \in V(\hbar),$
\begin{equation} \label{lower}
( Op_{\hbar}(\chi_{1}^{\delta}) \phi_{\mu}, \phi_{\mu} ) \geq
C(\epsilon) >0.
\end{equation}
\noindent Thus,
\begin{equation}
\| \phi_{\mu} \|^{2}_{L^{\infty}} \cdot \left( \int_{M}
\chi^{\delta}_{1} (x;\hbar ) \, \, d vol(x) \right) \geq C(\epsilon) ,
\end{equation}
\noindent uniformly for  $\hbar \in (0,\hbar_{0}(\epsilon)]$.
Since
$$ \int_{M} \chi_{1}^{\delta}(x;\hbar) \,\, d vol(x) = {\mathcal O}(\hbar^{\delta(n-\ell)})$$
\noindent with $\hbar^{-1} \in$ Spec $-\sqrt{\Delta}$, the lower
bound coming from (\ref{lower}) is:
$$\| \phi_{\lambda_{j} } \|_{L^{\infty}} \geq C(\epsilon)
\lambda_{j}^{\frac{n-\ell}{4} - \epsilon}.$$ \noindent The proof
of Lemma \ref{BLOWUP} (i) is complete. Lemma \ref{BLOWUP} (ii)
then follows by applying the H\"{o}lder inequality in the estimate
(\ref{lower}).
\end{proof}

\vspace{6mm}

\noindent {\bf Remark}:  From   our earlier results in the case of
the quantum Euler top [T1], it seems  likely that in fact $\|
\phi_{\lambda_{j} } \|_{\infty} \geq C
\lambda_{j}^{\frac{n-\ell}{8}}$ in the elliptic case and $\|
\phi_{\lambda_{j} } \|_{\infty} \geq C
\lambda_{j}^{\frac{n-\ell}{4}}  (\log \lambda_{j})^{-\alpha} $ for
an appropriate $\alpha>0$ in the hyperbolic case.

\end{document}